\documentclass[a4paper]{article}
\usepackage[english]{babel}
\usepackage{amsmath,amssymb,amsfonts,amsthm}
\usepackage[all]{xy}
\usepackage{ifpdf}
\ifpdf
  \usepackage{hyperref}
\else
  \usepackage[hypertex]{hyperref}
\fi

\newcommand{\dual}{*}
\newcommand{\rdual}{\circ}
\newcommand{\C}{\mathbb{C}}
\newcommand{\R}{\mathbb{R}}
\newcommand{\N}{\mathbb{N}}
\newcommand{\Z}{\mathbb{Z}}

\newcommand{\tc}{\,:\,}

\newcommand{\T}{\mathbb{T}}

\numberwithin{equation}{section}

\newtheorem{prp}{Proposition}[section]
\newtheorem{thm}[prp]{Theorem}
\newtheorem{cor}[prp]{Corollary}
\newtheorem{lem}[prp]{Lemma}

\theoremstyle{remark}
\newtheorem{rem}{Remark}[section]

\title{Generalized Uncertainty Inequalities}
\author{Alessio Martini\\(Scuola Normale Superiore, Pisa)\\{\tt\small a.martini@sns.it}}

\begin{document}
\maketitle


\section{Introduction}
The uncertainty principle, which is a fundamental feature of quantum mechanical systems, from a mathematical point of view can be considered as a ``meta-theorem'' in harmonic analysis, which can be summed up as: a nonzero function and its Fourier transform cannot both be sharply localized.

This qualitative statement has large varieties of quantitative formulations, extensions and generalizations (see \cite{MR1448337} for a survey). Here we are interested in generalizations of one of the most common quantitative restatements of the uncertainty principle, namely Heisenberg-Pauli-Weyl (HPW) inequality: for every $\alpha,\beta > 0$ there exists $C_{\alpha,\beta}$ such that
\[\|f\|_2 \leq C_{\alpha,\beta} \| |x|^\alpha f\|_2^\frac{\beta}{\alpha+\beta} \||\xi|^\beta \widehat f\|_2^\frac{\alpha}{\alpha+\beta}\]
for all $f \in L^2(\R^n)$. The inequality can also be rewritten as
\[\|f\|_2 \leq C_{\alpha,\beta} \| |x|^\alpha f\|_2^\frac{\beta}{\alpha+\beta} \|(-\Delta)^{\beta/2} f\|_2^\frac{\alpha}{\alpha+\beta}\]
and in this form it is possible to discuss its validity in more general contexts than $\R^n$ (e.g.\ in Riemannian manifolds, with $|x|$ interpreted as the distance from a fixed point and $\Delta$ as the Laplace-Beltrami operator).

The work~\cite{MR2355602} goes in this direction, obtaining uncertainty inequalities in ``spaces with polynomial volume growth'': measure spaces $(X,m)$ with a given ``distance-from-a-point'' function $\rho$ (which we can assume to be simply a non-negative measurable function on $X$) such that the measure of the ``balls'' (sublevel sets) $B_r = \{\rho < r\}$ is majorized by powers of the radius $r$:
\[m(B_r) \lesssim \begin{cases}
r^{q_0} &\text{for $r \leq 1$}\\
r^{q_\infty} &\text{for $r \geq 1$}
\end{cases}\]
for some $q_0,q_\infty \in \left]0,+\infty\right[$. In such a setting they obtain uncertainty inequalities of the form
\[\|f\|_2 \leq C_{\alpha,\beta} \| \rho^\alpha f\|_2^\frac{\beta}{\alpha+\beta} \|L^\beta f\|_2^\frac{\alpha}{\alpha+\beta}\]
where $L$ is any positive self-adjoint operator on $L^2(X,m)$ whose exponential semigroup $e^{-tL}$ satisfies the following ultracontractivity condition:
\begin{equation}\label{eq:ultracontractivity}
\|e^{-t L}\|_{1 \to \infty} \lesssim \begin{cases}
t^{-q_0} &\text{for $t \leq 1$}\\
t^{-q_\infty} &\text{for $t \geq 1$}.
\end{cases}
\end{equation}
Their proof gives also a ``local uncertainty inequality'' (from which the ``global'' one is derived):
\begin{equation}\label{eq:localhpwexp}
\|e^{-t L} f\|_2 \leq C_\alpha t^{-\alpha} \|\rho^\alpha f\|_2
\end{equation}
for $t$ small and $\alpha < q_0/2$, or for $t$ large and $\alpha < q/2$ (where $q = \min\{q_0,q_\infty\}$).

A first question which arises from this work is if the ``symmetry'' of the two factors in HPW inequality, given in $\R^n$ by Fourier transform, can be recovered, at least partially, in this more general setting.

Another question is if the ``polynomial growth'' condition can be relaxed, to include e.g.\ spaces with exponential volume growth, and what conditions must be satisfied in this case by the operator $L$.

The first problem is addressed specifically in \cite{MR2305079}, where a ``companion'' inequality of \eqref{eq:localhpwexp} is proved, i.e.
\begin{equation}\label{eq:companion}
\|e^{-r \rho} f\|_2 \leq C_\alpha r^{-\alpha} \|L^\alpha f\|_2
\end{equation}
for $r$ small and $\alpha < q_\infty/2$, or for $r$ large and $\alpha < q/2$. In this estimate the roles of the operator $L$ and the operator ``multiplication by $\rho$'' are swapped. The proof of \eqref{eq:companion} in \cite{MR2305079} is formally different from that of \eqref{eq:localhpwexp}, but the leading ideas are the same.

This suggests that the operator ``multiplication by $\rho$'' can be substituted with a generic positive self-adjoint operator $T$ on $L^2(X,m)$ (it should also be remarked that, by the spectral theorem, every self-adjoint operator can in fact be thought as a multiplication operator on some $L^2$ space). Let $F$ be the spectral measure associated to $T$ and set $F_r = F(\left[0,r\right[)$ for $r \geq 0$. Observing that, in the case of the multiplication operator $T f = \rho f$ we have also $F_r f = \chi_{B_r} f$ and $m(B_r) = \|\chi_{B_r}\|_1 = \|F_r\|_{\infty \to 1}$, the volume growth condition can be rewritten as
\begin{equation}\label{eq:grcnd1}
\|F_r\|_{\infty \to 1} \lesssim \begin{cases}
r^{q_0} &\text{for $r \leq 1$}\\
r^{q_\infty} &\text{for $r \geq 1$}
\end{cases}
\end{equation}
and in this form it makes sense also for a generic $T$.

To get now a similar condition on the operator $L$, we can use another remark in \cite{MR2305079}, the inequalities
\[e^{-1} \|E_t f\|_2 \leq \|e^{-L/t} f\|_2 \leq (e-1) \sum_{j=1}^\infty e^{-j} \|E_{jt} f\|_2\]
which imply the equivalence of condition \eqref{eq:ultracontractivity} to the following
\begin{equation}\label{eq:grcnd2}
\|E_t\|_{1 \to \infty} \lesssim \begin{cases}
t^{q_\infty} &\text{for $t \leq 1$}\\
t^{q_0} &\text{for $t \geq 1$}
\end{cases}
\end{equation}
It should also be noted that the thesis, the local uncertainty inequality (and its ``companion''), can be equivalently rephrased in terms of spectral projections:
\[\|E_{1/t} f\|_2 \leq C_\alpha t^{-\alpha} \|T^\alpha f\|_2\]
\[\|F_{1/r} f\|_2 \leq C_\alpha r^{-\alpha} \|L^\alpha f\|_2\]
while the global inequality takes the form
\[\|f\|_2 \leq C_{\alpha,\beta} \| T^\alpha f\|_2^\frac{\beta}{\alpha+\beta} \|L^\beta f\|_2^\frac{\alpha}{\alpha+\beta}\]
which is undoubtedly more ``symmetric''.

As they are now written, the uncertainty inequalities make sense not only in $L^2$ but also in a generic Hilbert space $H$. The problem is how to rephrase the growth hypotheses on spectral measures, since they are in terms of $L^1$ and $L^\infty$, which are Banach spaces having a close relationship with each other (duality) and with $L^2$. A suitable generalization is given by the concept of \emph{Banach couple} (see~\cite{MR0482275}): a pair $(X_0,X_1)$ of Banach spaces which are both (continuously) contained in a (Hausdorff) topological vector space $Z$ (so that we can also consider the sum $X_0+X_1$ and the intersection $X_0 \cap X_1$ as subspaces of $Z$). In fact, we will be interested in Banach couples which are \emph{regular} ($X_0 \cap X_1$ is dense in both $X_i$), \emph{reflexive} (in a sense which will be precised later) and with $X_0 = H$. For instance, if we choose $(L^2,L^p)$ as Banach couple (for $1 \leq p < \infty$), the growth hypotheses take the form of estimates on the norms $\|E_t\|_{p \to p'}$, $\|F_r\|_{p'\to p}$ (where $1/p+1/p' = 1$), so that the original case is recovered for $p = 1$. The case $p = \infty$ can be considered too, by \emph{regularization} of the couple $(L^2,L^\infty)$, i.e.\ by restricting to the couple $(L^2,L^\infty_0)$, where $L^\infty_0$ is the closure of $L^2 \cap L^\infty$ in $L^\infty$.

We now come to the second question, about the possibility of relaxing the growth conditions~\eqref{eq:grcnd1}, \eqref{eq:grcnd2} to include more general ``volume growths''. The first idea is that, as in the case of polynomial growth, the estimates on spectral projections of $L$ and $T$ should ``balance one another'', something like
\[\|E_{1/t}\|_{V \to V^\dual} \|F_{\eta t}\|_{V^\dual \to V} \lesssim 1\]
for some $\eta > 0$ and all $t$ (where $V = X_1$ in the Banach couple). In fact, what we require in the general case is that
\[\|F_r\|_{V^\dual \to V} \leq \Phi(r) \qquad \text{and} \qquad \|E_{1/t}\|_{V \to V^\dual} \Phi(\eta t) \lesssim 1\]
for some nonnegative measurable function $\Phi$ on $\left[0,+\infty\right[$ which satisfies the following admissibility hypothesis:
\[\int_0^r s^{-\gamma} \Phi(s) \,\frac{ds}{s} \lesssim r^{-\gamma} \Phi(r)\]
for some $\gamma > 0$ and all $r > 0$. This condition is satisfied by polynomial growth ($\Phi(r) = r^d$ with $d > \gamma$) but also by faster and slower growths (exponential, logarithmic).

In the following, local and global uncertainty inequalities are proved in this general context. The result is then applied to Riemannian manifolds (with Riemannian distance and Laplace-Beltrami operator), obtaining HPW inequalities on homogeneous simply connected manifolds with negative sectional curvature, on Riemannian symmetric spaces of non-compact type and, by restricting to the orthogonal complement of the kernel of the Laplacian, also on compact manifolds. Finally, similar results are obtained in the context of homogeneous graphs (with graph distance and difference Laplacian) and unimodular Lie groups (with Carnot-Carath\'eodory distances and left-invariant sublaplacians).

\subsection*{Acknowledgements}
I thank prof.\ Fulvio Ricci for his great expertise and helpfulness.
\section{Uncertainty inequalities}

\subsection{Preliminaries}\label{subsection:preliminaries}

From now on, all Banach spaces will be complex.

If $V$ is a Banach space, let $V^\dual$ denote the \emph{conjugate-dual} of $V$, i.e.\ the Banach space of continuous conjugate-linear functionals on $V$. If $F : V \to W$ is a continuous linear map of Banach spaces, let $F^\dual : W^\dual \to V^\dual$ be the \emph{transpose} of $F$, defined by $F^\dual(\phi) = \phi \circ F$. It is easy to see that $V^{\dual\dual}$ is naturally isomorphic to the linear bidual of $V$. Moreover, Riesz representation theorem for Hilbert spaces can be rephrased as follows: the map
\[H \ni v \mapsto \langle v, \cdot \rangle_H \in H^\dual\]
is a natural isometric linear isomorphism between any Hilbert space $H$ and its conjugate-dual (where naturality means that the transpose $F^\dual$ of a linear map $F$ between Hilbert spaces corresponds to the adjoint of $F$).

A \emph{Banach couple}\footnote{For a reference about Banach couples and Doolittle diagrams see \cite{MR0482275}, \cite{MR1107298}.} is a pair $(X_0,X_1)$ of Banach spaces which are both continuously included in a (Hausdorff) topological vector space $Z$; in this case, we can then form the intersection $X_0 \cap X_1$ and the sum $X_0 + X_1$ as subspaces of $Z$, which are also Banach spaces with suitable norms\footnote{A common choice is
\[\|x\|_{X_0 \cap X_1} = \max \{\|x\|_{X_0}, \|x\|_{X_1}\},\]
\[\|x\|_{X_0 + X_1} = \inf \{\|x_0\|_{X_0} + \|x_1\|_{X_1} \tc x_i \in X_i, \,x_0 + x_1 = x\}.\]}, so that the following diagram of inclusions
\[\xymatrix{
X_0 \cap X_1 \ar[r]^{i_0}\ar[d]^{i_1} & X_0 \ar[d]^{j_0}\\
X_1 \ar[r]^{j_1} & X_0+X_1
}\]
is both a pullback and a pushout (i.e.\ a so-called \emph{Doolittle diagram}).

A Banach couple $(X_0,X_1)$ is said \emph{regular} if $X_0 \cap X_1$ is dense in both $X_0,X_1$, or equivalently if both $X_0,X_1$ are dense in $X_0 + X_1$. In this case, all the maps in the conjugate-dual Doolittle diagram
\[\xymatrix{
(X_0 \cap X_1)^\dual \ar@{<-}[r]^{i_0^\dual}\ar@{<-}[d]^{i_1^\dual} & X_0^\dual \ar@{<-}[d]^{j_0}\\
X_1^\dual \ar@{<-}[r]^{j_1^\dual} & (X_0+X_1)^\dual
}\]
are injective, so that, by identifying $X_0^\dual,X_1^\dual$ with their images in $(X_0 \cap X_1)^\dual$, we can think of $(X_0^\dual,X_1^\dual)$ as a Banach couple, with $X_0^\dual + X_1^\dual = (X_0 \cap X_1)^\dual$, $X_0^\dual \cap X_1^\dual = (X_0 + X_1)^\dual$.

The conjugate-dual $(X_0^\dual,X_1^\dual)$ of a regular Banach couple $(X_0,X_1)$ need not be regular: $X_0^\dual \cap X_1^\dual$ is always weakly$^\dual$ dense in both $X_0^\dual,X_1^\dual$, but in general it is not strongly dense (however, if $X_i$ is reflexive, then $X_0^\dual \cap X_1^\dual$ is strongly dense in $X_i^\dual$). We can then consider the \emph{regularized conjugate-dual couple} $(X_0^\rdual, X_1^\rdual)$, where $X_i^\rdual$ is the closure in $X_i$ of $X_0 \cap X_1$.

Repeating this procedure, we obtain the regularized conjugate-bidual couple $(X_0^{\rdual\rdual}, X_1^{\rdual\rdual})$ and, as in the case of single Banach spaces, there are canonical continuous immersions $J_i : X_i \to X_i^{\rdual\rdual}$, defined by
\[J_i(x)(\phi) = \overline{\phi(x)},\]
which together are a morphism of Banach couples ($J_0|_{X_0 \cap X_1} = J_1|_{X_0 \cap X_1}$); if this morphism is an isomorphism (i.e.\ if both $J_i$ are isomorphisms) then the couple $(X_0,X_1)$ will be called \emph{reflexive}.

The notion of canonical immersion in the bidual for regular Banach couples is not perfectly analogous to the corresponding notion for single Banach spaces. The main differences are the following.

\begin{itemize}
\item In general the immersions $J_i : X_i \to X^{\rdual\rdual}_i$ are continuous and injective, but not necessarily isometric, nor homeomorphisms with their images. In fact, for $x \in X_i$, the norm of $J_i(x)$ in $X_i^{\rdual\rdual}$ is given by
\[p_i(x) = \sup_{0 \neq \phi \in X_i^\rdual} \frac{|\phi(x)|}{\|\phi\|_{X_i^\dual}}\]
which is a norm on $X_i$, since $X_i^\rdual$ is weakly$^\dual$ dense in $X_i^\dual$, but is not necessarily equivalent to the original norm $\|\cdot\|_{X_i}$. Since
\[p_i(x) = \sup_{t > 0} K_i(t,x) \qquad\text{for $x \in X_i$,}\]
where $K_i$ is the \emph{Peetre $K$-functional}
\[K_i(t,x) = \inf\{ \|x_i\|_{X_i} + t\|x_{1-i}\|_{X_{1-i}} \tc x_j \in X_j, \,x_0 + x_1 = x\}\]
(for $t > 0$, $x \in X_0 + X_1$), this inequivalence of norms occurs exactly when $X_i$ is not \emph{relatively complete} in $X_0 + X_1$, i.e.\ when the closed unit ball of $X_i$ is not closed in $X_0 + X_1$ (see~\S2.2 in \cite{MR1107298}).

\item If both $X_i$ are reflexive, then the couple $(X_0,X_1)$ is reflexive too. If one of the $X_i$ is reflexive, then $(X_0,X_1)$ need not be reflexive, but $(X_0^\rdual,X_1^\rdual)$ is certainly reflexive.
\end{itemize}

In the following, we will in fact be interested in reflexive regular Banach couples of the form $(H,V)$, where $H$ is a Hilbert space. In this case, modulo identification by Riesz representation theorem, $(H,V^\rdual)$ is the regularized conjugate-dual couple; moreover, by replacing the norm of $V$ with the equivalent norm on $V^{\rdual\rdual}$, we can always suppose that the immersion $V \to V^{\rdual\rdual}$ is an isometry, so that we can identify $V^{\rdual\rdual}$ with $V$. Under these hypotheses, we have:

\begin{lem}\label{lem:adjnorm}
Let $P : H \to H$ be a continuous linear operator. The following are equivalent:
\begin{itemize}
\item $P$ is continuous $V \to H$,
\item $P^\dual$ is continuous $H \to V^\rdual$,
\item $P^\dual P$ is continuous $V \to V^\rdual$;
\end{itemize}
moreover
\[\|P\|^2_{V \to H} = \|P^\dual\|^2_{H \to V^\rdual} = \|P^\dual P\|_{V \to V^\rdual}.\]
\end{lem}
\begin{proof}
Let $P$ be continuous $V \to H$. Then, for every $v \in H$, $w \in H \cap V$,
\[|\langle w, P^\dual(v)\rangle| = |\langle P(w), v\rangle| \leq \|P(w)\|_H \|v\|_H \leq \|P\|_{V \to H} \|w\|_{V} \|v\|_H,\]
so that $P^\dual(v) \in H \cap V^\dual = H \cap V^\rdual$ and
\[\|P^\dual(v)\|_{V^\rdual} \leq \|P\|_{V \to H} \|v\|_H,\]
therefore $P^\dual$ is continuous $H \to V^\rdual$ and
\[\|P^\dual\|_{H \to V^\rdual} \leq \|P\|_{V \to H}.\]

Let now $P^*$ be continuous $H \to V^\rdual$, $v \in H \cap V$, $w \in H$. Then
\[|\langle P(v), w \rangle| = |\langle v, P^\dual(w)\rangle| \leq \|v\|_{V} \|P^\dual(w)\|_{V^\rdual} \leq \|P^\dual\|_{H \to V^\rdual} \|v\|_{V} \|w\|_H,\]
so that
\[\|P(v)\|_{H} \leq \|P^\dual\|_{H \to V^\rdual} \|v\|_{V},\]
which means that $P$ is continuous $V \to H$ and $\|P\|_{V \to H} \leq \|P^\dual\|_{H \to V^\rdual}$. Moreover, $P^\dual P$ is continuous $V \to V^\dual$ and
\[\|P^\dual P\|_{V \to V^\rdual} \leq \|P\|_{V \to H} \|P^\dual\|_{H \to V^\rdual} \leq \|P\|_{V \to H}^2.\]

Finally, let $P^\dual P$ be continuous $V \to V^\rdual$. For every $v \in H \cap V$,
\[\|P(v)\|_H^2 = \langle v,P^\dual P(v) \rangle \leq \|v\|_{V} \|P^\dual P(v)\|_{V^\rdual} \leq \|P^\dual P\|_{V \to V^\rdual} \|v\|_{V}^2,\]
so that $P$ is continuous $V \to H$ and $\|P\|_{V \to H}^2 \leq \|P^\dual P\|_{V \to V^\rdual}$.
\end{proof}

\subsection{The main theorems}\label{subsection:theorems}
If $H$ is a Hilbert space, $T$ is an unbounded self-adjoint operator on $H$, for every $f \in H$ which does not belong to the domain of $T$, we set $\|Tf\|_H = +\infty$, so that the equality
\[\|Tf\|_H = \sqrt{\int_\R \lambda^2 \,\|F(d\lambda) f\|_H^2}\]
(where $F$ is the spectral measure associated to $T$) holds for all $f \in H$.

In the following $(H,V)$ will be a reflexive regular Banach couple, where $H$ is a Hilbert space. Moreover, $L,T$ will denote (possibly unbounded) positive self-adjoint operators on $H$, $E,F$ the associated spectral measures and, for all $\lambda \geq 0$,
\[E_\lambda = E(\left[0,\lambda\right[), \qquad F_\lambda = F(\left[0,\lambda\right[).\]

\begin{thm}\label{thm:uncertainty1}
Let $A = \left]a,b\right[ \subseteq \left]0,+\infty\right[$ an open interval, $\eta,\delta > 0$. Let $\Phi$ a nonnegative measurable function on $B = [0,\eta b^\delta[$ such that
\begin{equation}\label{eq:ip0}
\|F_r\|_{V^\rdual \to V} \leq \Phi(r) \qquad\text{for all $r \in B$}
\end{equation}
and that, for some $K > 0$,
\begin{equation}\label{eq:ip2}
\|E_{1/t}\|_{V \to V^\rdual} \,\,\Phi(\eta t^\delta) \leq K^2 \qquad\text{for all $t \in A$.}
\end{equation}

Moreover, suppose that, for some $\gamma,M > 0$,
\begin{equation}\label{eq:ip1}
\int_0^r s^{-2\gamma} \Phi(s) \,\frac{ds}{s} \leq M r^{-2\gamma} \,\Phi(r) \qquad\text{for all $r \in \eta A^\delta$,}
\end{equation}
where $\eta A^\delta = \left]\eta a^\delta,\eta b^\delta\right[$.

Then, for all $f \in H$,
\[\|E_{1/t} f\|_H \leq C t^{-\gamma\delta} \|T^\gamma f\|_H \qquad\text{for all $t \in A$,}\]
where $C = \eta^{-\gamma}(1 + K \sqrt{1+2\gamma M})$.
\end{thm}
\begin{proof}
Let $f \in H$ be in the domain of $T^\gamma$, $t,r > 0$ and set $f_r = F_r f$, $f = f_r + f^r$, so that
\[\|E_{1/t} f\|_H \leq \|E_{1/t} f^r\|_H + \|E_{1/t} f_r\|_H.\]

We have immediately
\[\|E_{1/t} f^r\|_H \leq \|f^r\|_H \leq \|(1-F_r)T^{-\gamma}\|_{H\to H} \|T^\gamma f\|_H \leq r^{-\gamma} \|T^\gamma f\|_H.\]

Note that, for $g \in H \cap V^\dual$, if $\nu = \|F(\cdot) g\|_H^2$, for every $x \in B$
\[\nu(\left[0,x\right[) = \langle g, F_x g \rangle \leq \|F_x\|_{V^\rdual \to V} \|g\|_{V^\dual}^2 \leq \|g\|_{V^\dual}^2 \Phi(x)\]
by \eqref{eq:ip0}, therefore, if $r \in \eta A^\delta$, integrating by parts,
\begin{multline*}
\|F_r T^{-\gamma} g\|_H^2 = \int_{\left[0,r\right[} s^{-2\gamma} \,\nu(ds) = r^{-2\gamma} \nu(\left[0,r\right[) + 2\gamma\int_0^r s^{-2\gamma} \nu(\left[0,s\right[) \,\frac{ds}{s} \\
\leq \|g\|_{V^\dual}^2 \Phi(r) + 2\gamma \|g\|_{V^\dual}^2 \int_0^r s^{-2\gamma} \Phi(s) \,\frac{ds}{s} \leq (1+2\gamma M) r^{-2\gamma} \Phi(r) \|g\|_{V^\dual}^2
\end{multline*}
by \eqref{eq:ip1}.

Since $T^\gamma f$ is in the domain of $T^{-\gamma}$ and $F_r(H) \subseteq V$ by \eqref{eq:ip0} and Lemma~\ref{lem:adjnorm}, then $f_r = F_r T^{-\gamma} T^\gamma f \in V$; moreover, for every $g \in H \cap V^\dual$,
\begin{multline*}
|\langle g, f_r\rangle| = |\langle F_r T^{-\gamma} g,T^\gamma f\rangle| \leq \|F_r T^{-\gamma} g\|_H \|T^\gamma f\|_H \\
\leq M' r^{-\gamma} \sqrt{\Phi(r)} \|g\|_{V^\dual} \|T^\gamma f\|_H,
\end{multline*}
where $M' = \sqrt{1+2\gamma M}$, so
\[\|f_r\|_V \leq M' r^{-\gamma} \sqrt{\Phi(r)} \|T^\gamma f\|_H,\]
therefore
\[\|E_{1/t} f_r\|_H \leq \|E_{1/t}\|_{V \to H} \|f_r\|_V \leq r^{-\gamma} M' \sqrt{\left\|E_{1/t}\right\|_{V \to V^\rdual} \Phi(r)} \|T^\gamma f\|_H\]
by Lemma~\ref{lem:adjnorm}.

Putting all together,
\[\|E_{1/t} f\|_H \leq r^{-\gamma} \left(1 + M' \sqrt{\left\|E_{1/t}\right\|_{V \to V^\rdual} \Phi(r)}\right) \|T^\gamma f\|_H,\]
so that, choosing $r = \eta t^\delta$, $t \in A$, we get the result by \eqref{eq:ip2}.
\end{proof}

\begin{rem}\label{rem:expproj}
The inequalities
\[e^{-1} \chi_{\left[0,t\right[}(\lambda) \leq e^{-\lambda/t} \leq (e-1) \sum_{j=1}^\infty e^{-j} \chi_{\left[0,jt\right[}(\lambda),\]
true for all $\lambda \geq 0$, imply that, for all $f \in H$,
\[e^{-1} \|E_t f\|_H \leq \|e^{-L/t} f\|_H \leq (e-1) \sum_{j=1}^\infty e^{-j} \|E_{jt} f\|_H.\]
In particular, by Lemma~\ref{lem:adjnorm},
\[\|E_t\|_{V \to V^\rdual} = \|E_t\|^2_{V \to H} \leq e^2 \|e^{-L/t}\|^2_{V \to H} = e^2 \|e^{-2L/t}\|^2_{V \to V^\rdual}\]
and, analogously,
\[\|F_t\|_{V^\rdual \to V} \leq e^2 \|e^{-2T/t}\|_{V^\rdual \to V},\]
so that, in the hypotheses of the previous theorem, the estimates of the operator norms of the spectral measures $E,F$ can be replaced\footnote{Note that, in case of non-polynomial growth, the hypotheses on the spectral measures are weaker than the corresponding hypotheses on the semigroups. Moreover, estimates on spectral measures are easier to be managed when the operator is somehow rescaled.} by analogous estimates of the norms of the semigroups generated by $L,T$. Moreover, also the thesis can be rewritten in terms of the semigroup generated by $L$, because from
\[\|E_{1/t} f\|_H \leq C t^{-\gamma\delta} \|T^\gamma f\|_H\]
it follows that
\[\|e^{-tL} f\|_H \leq C' t^{-\gamma\delta} \|T^\gamma f\|_H,\]
where
\[C' = C (e-1) \sum_{j=1}^\infty e^{-j} j^{\gamma\delta} < +\infty.\]
\end{rem}

\begin{rem}\label{rem:afortiori}
If \eqref{eq:ip1} holds for some $\gamma > 0$, then it holds also for every $\gamma' < \gamma$, since
\begin{multline*}
\int_0^r s^{-2\gamma'} \Phi(s)\,\frac{ds}{s} \leq r^{2(\gamma-\gamma')} \int_0^r s^{-2\gamma} \Phi(s)\,\frac{ds}{s}\\
\leq r^{2(\gamma-\gamma')} M r^{-2\gamma} \Phi(r) = M r^{-2\gamma'} \Phi(r).
\end{multline*}
\end{rem}

\begin{thm}\label{thm:uncertainty2}
Suppose that, for some $f \in H$, $\gamma,\delta,C > 0$,
\[\|E_{1/t} f\|_H \leq C t^{-\gamma \delta} \|T^\gamma f\|_H \qquad\text{for all $t > 0$.}\]
Then, for all $\alpha \geq \gamma$, $\beta > 0$,
\begin{equation}\label{eq:hpw}
\|f\|_H \leq D_{\alpha,\beta} \|T^\alpha f\|_H^{\frac{\beta}{\alpha+\beta}} \|L^{\beta\delta} f\|_H^{\frac{\alpha}{\alpha+\beta}},
\end{equation}
where $D_{\alpha,\beta} > 0$ depends only on $C,\gamma,\alpha,\beta$.
\end{thm}
\begin{proof}
Suppose first $\alpha = \gamma$. Then, for all $t > 0$, by the spectral theorem
\begin{multline*}
\|f\|_H \leq \|E_{1/t} f\|_H + \|(1-E_{1/t}) L^{-\beta\delta} L^{\beta\delta} f\|_H \\
\leq (1+C) \left(t^{-\gamma\delta}\|T^\gamma f\|_H + t^{\beta\delta} \|L^{\beta\delta} f\|_H \right),
\end{multline*}
from which, optimizing in $t$, we obtain \eqref{eq:hpw} with $D_{\gamma,\beta} = (1+C) (\gamma/\beta)^{\frac{\beta-\gamma}{\gamma+\beta}}$.

Let now $\alpha > \gamma$. Then, for all $f \in H$, $\epsilon > 0$, if $\nu = \|F(\cdot)f\|^2_H$,
\begin{multline*}
\epsilon^{-\gamma} \|T^{\gamma} f\|_H = \sqrt{\int_{\R^+} (\lambda/\epsilon)^{2\gamma} \,d\nu(\lambda)} \leq \sqrt{\int_{\R^+} (1+(\lambda/\epsilon)^\alpha)^2 \,d\nu(\lambda)}\\
= \|(1+\epsilon^{-\alpha}T^\alpha)f\|_H \leq \|f\|_H + \epsilon^{-\alpha}\|T^{\alpha} f\|_H,
\end{multline*}
hence, optimizing in $\epsilon$,
\[\|T^\gamma f\|_H \leq K_{\alpha,\gamma} \|f\|_H^{1-\frac{\gamma}{\alpha}} \|T^\alpha f\|_H^{\frac{\gamma}{\alpha}}\]
(where $K_{\alpha,\gamma} = (\alpha/\gamma-1)^\frac{2\gamma-\alpha}{\alpha}$). Plugging this into \eqref{eq:hpw} with $\alpha$ replaced by $\gamma$, we obtain
\[\|f\|^{\gamma+\beta}_H \leq D^{\gamma+\beta}_{\gamma,\beta} K^{\beta}_{\alpha,\gamma} \|f\|_H^{\left(1-\frac{\gamma}{\alpha}\right)\beta} \|T^{\alpha} f\|_H^{\frac{\gamma}{\alpha}\beta} \|L^{\beta\delta} f\|_H^\gamma,\]
that is \eqref{eq:hpw} with $D_{\alpha,\beta} = D_{\gamma,\beta}^{\frac{\alpha}{\gamma} \frac{\gamma+\beta}{\alpha+\beta}} K_{\alpha,\gamma}^{\frac{\alpha}{\gamma} \frac{\beta}{\alpha+\beta}}$.
\end{proof}

As it is formulated, Theorem~\ref{thm:uncertainty2} shows that global uncertainty inequalities can be obtained directly from local ones, which must hold for all times $t > 0$ but can be limited only to a certain subset of $H$. This formulation can be useful when local uncertainty inequalities are obtained by other means than Theorem~\ref{thm:uncertainty1}. However, we can certainly put together Theorems~\ref{thm:uncertainty1}, \ref{thm:uncertainty2} obtaining

\begin{cor}
In the hypotheses of Theorem~\ref{thm:uncertainty1} with $A = \left]0,+\infty\right[$, for all $\alpha,\beta > 0$, $f \in H$,
\[\|f\|_H \leq D_{\alpha,\beta} \|T^\alpha f\|_H^{\frac{\beta}{\alpha+\beta}} \|L^{\beta\delta} f\|_H^{\frac{\alpha}{\alpha+\beta}},\]
where $D_{\alpha,\beta} > 0$ depends only on $M,K,\eta,\gamma,\alpha,\beta$.
\end{cor}

\subsection{The hypothesis on the growth}\label{subsection:discussion}

The importance of \eqref{eq:ip1} is in that it allows to separate in two distinct factors the dependence on $\Phi$ and the dependence on $\gamma$ (so that hypothesis \eqref{eq:ip2} does not depend on $\gamma$).

In order to simplify the form of the hypothesis, we set $\alpha = 2\gamma$, $I = \eta A^\delta$, $C_{I,\alpha} = M$. The inequality then becomes
\begin{equation}\label{eq:ip1bis}
\int_0^r s^{-\alpha} \Phi(s) \,\frac{ds}{s} \leq C_{I,\alpha} \, r^{-\alpha} \Phi(r) \qquad\text{for all $r \in I$.}
\end{equation}
We are now going to discuss necessary or sufficient conditions for the existence of $C_{I,\alpha} > 0$ such that \eqref{eq:ip1bis} holds, where $\alpha > 0$, $I \subseteq \left]0,+\infty\right[$ is a non-empty interval and $\Phi$ a finite non-null non-negative measurable function defined on an interval $B \subseteq \left[0,+\infty\right[$ containing $I \cup \{0\}$.

In remark~\ref{rem:afortiori} we have already pointed out that, if \eqref{eq:ip1bis} holds for some $\alpha > 0$, it holds also for all $\alpha' > 0$ smaller than $\alpha$ with $C_{I,\alpha'} = C_{I,\alpha}$.

First of all, since $\Phi$ is finite, a necessary condition for \eqref{eq:ip1bis} to hold is that
\begin{equation}\label{eq:l1loc}
\int_0^\epsilon \frac{\Phi(s)}{s^{\alpha+1}} \,ds < +\infty \qquad\text{for some $\epsilon > 0$.}
\end{equation}

If $\sup I = +\infty$, information on the behavior of $\Phi$ in a neighborhood of $+\infty$ can also be recovered. In fact, since $\Phi \neq 0$, $\Phi \geq 0$, there exists $r' \in I$ such that
\[C_{I,\alpha}\, r^{-\alpha} \Phi(r) \geq \int_0^{r'} s^{-\alpha} \Phi(s) \,\frac{ds}{s} > 0 \qquad\text{for all $r \geq r'$,}\]
by \eqref{eq:ip1bis}, i.e.
\[\Phi(r) \gtrsim r^{\alpha} \qquad\text{for $r \to +\infty$.}\]

Suppose now that \eqref{eq:l1loc} holds and moreover that $\Phi$ is absolutely continuous, so that it admits a distributional derivative $\Phi' = f$ which is $L^1_{\mathrm{loc}}(B)$. In this case, \eqref{eq:ip1bis} becomes
\[\int_0^r f(s)s^{-\alpha}\,ds \leq C'_{I,\alpha} \, r^{-\alpha} \int_0^r f(s) \,ds\]
(where $C'_{I,\alpha} = 1 + 2\alpha C_{I,\alpha}$).

If $f(s) = s^{d-1}$ for some $d > 0$ and for $s$ small, then \eqref{eq:ip1bis} holds for $r\to 0^+$ iff $\alpha < d$.

If $f(s) = s^{d-1}$ for some $d > 0$ and for $s$ large, then \eqref{eq:ip1bis} holds for $r \to +\infty$ iff $\alpha < d$.

Another sufficient condition for \eqref{eq:ip1bis} to hold for $r\to +\infty$ is that $f(s)s^{-\alpha}$ is definitely nondecreasing; in fact, if $f(s) s^{-\alpha}$ is nondecreasing for $s > r_0 \geq 0$, then
\[\int_{r_0}^r f(s) s^{-\alpha} \,ds = \int_{r_0}^{r/2} + \int_{r/2}^r \leq 2 \int_{r/2}^r \leq 2^{\alpha+1} r^{-\alpha} \int_0^r f(s)\,ds\]
for all $r > 2r_0$. Moreover, note that, if $f(s) s^{-\alpha}$ is nondecreasing in a neighborhood of $0$, the same argument proves \eqref{eq:ip1bis} for $r\to 0^+$.

A case not included in the previous ones in which \eqref{eq:ip1bis} still holds for $r\to +\infty$ is $f(s) = (\log s)^\delta$ for $s$ large, $\delta > 0$, $0 < \alpha < 1$, since integrating by parts it is easily obtained that
\[\int_1^r s^{-\alpha} (\log s)^\delta \,ds \lesssim r^{1-\alpha} (\log r)^\delta \asymp r^{-\alpha} \int_1^r (\log s)^\delta \,ds \qquad\text{for $r \to +\infty$.}\]

%
%
%

\subsection{Hilbert-Banach couples of Lebesgue spaces}
From what we said in \S\ref{subsection:preliminaries}, it is clear that the hypotheses of \S\ref{subsection:theorems} on the regular Banach couple $(H,V)$ are satisfied if $H$ is Hilbert and $V$ is reflexive (and in this case $V^\rdual = V^\dual$). In particular, fixed a measure space $(X,m)$, those hypotheses are certainly satisfied by the couple of Lebesgue spaces $(L^2,L^p)$ on $(X,m)$ for $1 < p < \infty$.

Let us consider now the case $p = 1$, that is the couple $(L^2,L^1)$. This is certainly a regular Banach couple. Moreover, if
\[L^\infty_\sigma = \{f \in L^\infty \tc \text{$f$ is null out of a $\sigma$-finite subset of $X$}\}\]
(we are not supposing that $m$ is $\sigma$-finite), then $L^\infty_\sigma$ is a closed subspace of $L^\infty$, $(L^1)^\dual$ contains isometrically $L^\infty_\sigma$ as a subspace and $L^2 \cap (L^1)^\dual \subseteq L^\infty_\sigma$. Let $L^\infty_0$ be the closure in $(L^1)^\dual$ of this intersection, which is the closure in $L^\infty$ of the space of simple measurable functions of $(X,m)$ which are null out of a set of finite measure. Then $(L^2,L^\infty_0)$ is the regularized conjugate-dual of $(L^2,L^1)$.

Now, it is easy to see that $L^1$ is isometrically embedded in $(L^\infty_0)^\dual$ (since every $f \in L^1$ is null out of a $\sigma$-finite subset of $X$) and that $L^2 \cap (L^\infty_0)^\dual \subseteq L^1$; on the other hand, $L^2 \cap L^1$ is dense in $L^1$, therefore $L^1$ is the closure of $L^2 \cap (L^\infty_0)^\dual$ in $(L^\infty_0)^\dual$, so that $(L^2,L^1)$ is the regularized conjugate-dual of $(L^2,L^\infty_0)$. By a careful examination of the implicit identifications, it is then not difficult to see that $(L^2,L^1)$ is reflexive.

We have thus obtained that $(L^2,L^1)$, $(L^2,L^\infty_0)$ are both reflexive regular Banach couples $(H,V)$. Moreover, $V^\rdual = L^\infty_0$ in the former case, whereas in the latter $V^\rdual = L^1$. This shows an interesting mutual duality between $L^1$ and $L^\infty_0$, which holds in spite of non-reflexivity of the single Banach spaces and without any hypotheses of $\sigma$-finiteness of the measure.

\section{Applications}

\subsection{Uncertainty inequalities on Riemannian manifolds}\label{subsection:riemannian}

As we said in the introduction, Riemannian manifolds are a suitable setting to generalize uncertainty inequalities, since the notions of ``Laplacian'' and ``distance from a given point'' are meaningful there.

Let $M$ be a (connected) Riemannian manifold, $d$ the Riemannian metric, $m$ the Riemannian measure, $\Delta$ the Laplace-Beltrami operator. Chosen a point $x_0 \in M$, let $\rho = d(x_0,\cdot)$ and let $T$ be the operator ``multiplication by $\rho$''. Then $T$ is a positive self-adjoint operator on $L^2(M)$ and
\[\|F_r\|_{\infty \to 1} = \|\chi_{\{\rho < r\}}\|_1 = m(B(x_0,r)).\]

Suppose now that $M$ is a \emph{complete} Riemannian $n$-manifold. Then $L = -\Delta$, as an operator on $L^2(M)$, is (essentially) self-adjoint and positive (see~\cite{MR705991}); moreover the semigroup $e^{-tL}$ ($t > 0$) admits a kernel function $h_t$, the so-called \emph{heat kernel}, such that
\begin{itemize}
\item $(t,x,y) \mapsto h_t(x,y)$ is $C^\infty$ on $\left]0,+\infty\right[ \times M \times M$;
\item $h_t(x,y) > 0$ and $h_t(x,y) = h_t(y,x)$;
\item $h_t(x,y) \leq \sqrt{h_t(x,x) \,h_t(y,y)}$; 
\item $e^{-t L} f(x) = \int_M h_t(x,y) f(y) \,dm(y)$ for $m$-a.e. $x$.
\end{itemize}
In particular (cf.\ remark~\ref{rem:expproj})
\[\|E_{1/t}\|_{1 \to \infty} \lesssim \|e^{-2t L}\|_{1 \to \infty} = \|h_{2t}\|_\infty.\]

It is then interesting to see if the quantities $m(B(x_0,r))$ and $\|h_t\|_\infty$ are related in some way. In fact, there are several results (see e.g.\ \cite{MR1736868}) about the validity of the estimate
\begin{equation}\label{eq:lapvolestimate}
h_t(x,x) \,m(B(x,\sqrt{t})) \lesssim 1.
\end{equation}

First of all, \eqref{eq:lapvolestimate} always holds for small times $t>0$ locally in $x \in M$. This means that, if $M$ is e.g.\ compact or homogeneous, then \eqref{eq:lapvolestimate} holds uniformly on $M$ for small times. In this hypothesis, since $m(B(x_0,r)) \asymp r^n$ for $r \to 0^+$, it is sufficient to put $\Phi(r) = cr^n$ for a suitable $c > 0$ to get
\[\|F_r\|_{\infty \to 1} \leq \Phi(r),\qquad \|E_{1/t}\|_{1 \to \infty} \Phi(t^{1/2}) \lesssim 1 \qquad\text{for $r,t$ small},\]
and analogously, choosing $\Phi(r) = cr^{n/2}$, we get
\[\|E_r\|_{\infty \to 1} \leq \Phi(r),\qquad \|F_{1/t}\|_{1 \to \infty} \Phi(t^2) \lesssim 1 \qquad\text{for $r,t$ large}.\]
Therefore, by Theorem~\ref{thm:uncertainty1} and \S\ref{subsection:discussion} we obtain local uncertainty inequalities for small times: for $0 < \gamma < n/2$, $f \in L^2(M)$,
\begin{equation}\label{eq:local1}
\left.\begin{array}{r}
\|E_{1/t} f\|_2\\
\|e^{t\Delta} f\|_2
\end{array}\right\} \leq C_\gamma t^{-\gamma/2} \|\rho^\gamma f\|_2 \qquad \text{for $t$ small;}
\end{equation}
\begin{equation*}\label{eq:local2}
\left.\begin{array}{r}
\|\chi_{\{\rho < t\}} f\|_2\\
\|e^{-\rho/t} f\|_2
\end{array}\right\} \leq C_\gamma t^{\gamma} \|(-\Delta)^{\gamma/2} f\|_2 \qquad \text{for $t$ small.}
\end{equation*}

To get global uncertainty inequalities, in order to apply Theorem~\ref{thm:uncertainty2} we need to extend at least one of the local inequalities also to large times. If \eqref{eq:lapvolestimate} (or something similar) holds uniformly and for all times (see \cite{MR1736868} for sufficient conditions), if the rate of growth of the measure of the balls is independent of the center and moreover satisfies \eqref{eq:ip1}, then we can apply Theorem~\ref{thm:uncertainty1} also for large times.

A particularly simple case to be considered is when the Laplacian has a \emph{spectral gap}, i.e. the spectrum of $L$ is bounded from below by a constant $b > 0$. This holds e.g.\ when $M$ is simply connected and all sectional curvatures are bounded from above by a negative constant, by a result of McKean (see~\cite{MR705991}). In this cases, local inequalities for large times,
\begin{equation*}\label{eq:local3}
\left.\begin{array}{r}
\|\chi_{\{\rho < t\}} f\|_2\\
\|e^{-\rho/t} f\|_2
\end{array}\right\} \leq C_{\gamma,\delta} \,t^{\delta} \|(-\Delta)^{\gamma} f\|_2 \qquad \text{for $t$ large}
\end{equation*}
\begin{equation*}\label{eq:local4}
\left.\begin{array}{r}
\|E_{1/t} f\|_2\\
\|e^{t\Delta} f\|_2
\end{array}\right\} \leq C_{\gamma,\delta} \,t^{-\delta} \|\rho^\gamma f\|_2 \qquad \text{for $t$ large}
\end{equation*}
for all $\gamma,\delta > 0$, are trivially true (the former because $(-\Delta)^\gamma$ has a bounded inverse, the latter since $E_{1/t} = 0$ for $t$ large). Putting together the results for $t$ small and $t$ large and applying Theorem~\ref{thm:uncertainty2}, we obtain the following result, perfectly analogous to the Euclidean case:

\begin{cor}
If the Laplacian on the Riemannian manifold $M$ has a spectral gap, then, for all $\alpha,\beta > 0$, $f \in L^2(M)$,
\[\|f\|_2 \leq C_{\alpha,\beta} \|\rho^\alpha f\|_2^\frac{\beta}{\alpha+\beta} \|(-\Delta)^{\beta/2} f\|_2^\frac{\alpha}{\alpha+\beta},\]
\end{cor}

A different way to deal with a spectral gap is to replace $L$ with the operator $\tilde L = L - b$. In order to obtain results in this case we need precise information on the behavior of the norms of spectral projections $E_t$ of $L$ in a neighborhood of $b$, or at least on the decay of the heat kernel. Let us consider, for instance, a Riemannian symmetric spaces of non-compact type $M$ of dimension $n$ and rank $k$; chosen a system of positive roots, let $l$ be the norm of the sum of positive roots, counted with multiplicities, $s$ be the number of positive indivisible roots. Then it is known (see~\cite{MR2018351}) that $b = l^2/4 > 0$,
\[\|h_t\|_\infty \asymp t^{-\frac{n}{2}} e^{-bt} (1+t)^{\frac{n-k-2s}{2}} \asymp \begin{cases}
t^{-\frac{n}{2}} &\text{for $t \to 0^+$}\\
t^{-\frac{k+2s}{2}} e^{-bt} &\text{for $t \to +\infty$}
\end{cases}\]
whereas (cf.~\cite{MR1465601}, Theorem~6.2)
\[m(B(x_0,r)) \asymp \begin{cases}
r^n &\text{for $t \to 0^+$}\\
r^{\frac{k-1}{2}}e^{lr} &\text{for $t \to +\infty$}
\end{cases}\]
Since $e^{-t \tilde L} = e^{bt} e^{-t L}$, we have in particular
\[\|e^{-t \tilde L}\|_{1 \to \infty} \asymp \begin{cases}
t^{-\frac{n}{2}} &\text{for $t \to 0^+$}\\
t^{-\frac{k+2s}{2}} &\text{for $t \to +\infty$}
\end{cases}\]
To obtain uncertainty inequalities for $\tilde L$ we can then replace the distance function $\rho$ with
\[\tilde \rho = (1+\rho)^{\frac{k-1}{2(k+2s)}} e^{\frac{l}{k+2s} \rho} - 1,\]
so that
\[m(\{\tilde \rho < r\}) \lesssim \begin{cases}
r^n &\text{for $r \to 0^+$}\\
r^{k+2s} &\text{for $t \to +\infty$}
\end{cases}\]
Therefore local inequalities for all times and then global inequalities can be obtained for $\tilde\rho, \tilde L$ by applying Theorems~\ref{thm:uncertainty1}, \ref{thm:uncertainty2}:

\begin{cor}
If $M$ is a Riemannian symmetric space of non-compact type, for all $\alpha,\beta > 0$, $f \in L^2(M)$,
\[\|f\|_2 \leq C_{\alpha,\beta} \|((1+\rho)^\frac{k-1}{2(k+2s)} e^{\frac{l}{k+2s} \rho} -1)^\alpha f\|_2^\frac{\beta}{\alpha+\beta} \|(L-b)^{\beta/2} f\|_2^\frac{\alpha}{\alpha+\beta}.\]
\end{cor}

Note that, instead of ``exponentiating'' the distance function $\rho$, we could have ``taken the logarithm'' of the Laplacian $\tilde L$, thus getting another set of inequalities.

Another particular case is when $M$ is compact. Here, local inequalities for $\rho,L$ cannot be extended to large times, and global inequalities cannot hold, since the Laplacian has a non-null kernel, the space of constant functions on $M$ (which are in $L^2(M)$ if $M$ is compact). However, we can restrict to the orthogonal complement $H_0$ of $\ker L$, i.e.\ the space of functions with null mean value. Since $M$ is compact, the spectrum of $L$ is discrete (see~\cite{MR1736868}), so that $E_{1/t}|_{H_0} = 0$ for $1/t$ smaller than the first positive eigenvalue of $M$. Therefore \eqref{eq:local1} holds also for $t$ large if $f \in H_0$; then, by Theorem~\ref{thm:uncertainty2} we obtain:

\begin{cor}
If $M$ is a compact Riemannian manifold, for all $\alpha,\beta > 0$, $f \in L^2(M)$ with null mean value,
\[\|f\|_2 \leq C_{\alpha,\beta} \|\rho^\alpha f\|_2^\frac{\beta}{\alpha+\beta} \|(-\Delta)^{\beta/2} f\|_2^\frac{\alpha}{\alpha+\beta}.\]
\end{cor}

\subsection{Uncertainty inequalities on graphs}

A considerably studied subject is the spectral theory of graphs (see e.g.\ \cite{MR986363} for a survey). On a (unoriented multi)graph $G = (V,E)$ there are a canonical distance function $d$ on vertices (given by the minimum length of a path joining two vertices), a canonical measure $m$ (the counting measure, which is a Borel measure with respect to the discrete topology induced by $d$ on $V$) and, if $G$ is locally finite (i.e.\ $\deg(u) < \infty$ for all vertices $u$, where $\deg(u)$ is the number of edges emanating from $u$), two \emph{difference Laplacians}:
\[\Delta_A = D - A, \qquad \Delta_P = I - P,\]
where $A$ is the \emph{adjacency matrix} of $G$ (i.e.\ $a_{uv}$ is the number of edges between $u$ and $v$), $D = (\delta_{uv} \deg(u))_{u,v \in V}$, $P = (a_{uv}/\deg (u))_{u,v \in V}$ is the \emph{transition matrix} of $G$ and $I = (\delta_{uv})_{u,v \in V}$ is the identity matrix.

Supposing $G$ homogeneous (i.e.\ $\deg(u)$ is independent of $u$ and denoted by $\deg(G)$) and locally finite, then
\[\Delta_A = \deg(G) \Delta_P, \qquad D = \deg(G) I,\]
so that $A$ is a bounded self-adjoint operator on $L^2(G)$, with norm at most $\deg (G)$, and spectral information on $A$ carries over to $\Delta_A, P, \Delta_P$.

In these hypotheses, let $x_0 \in V$, $\rho = d(x_0,\cdot)$, $T$ the operator ``multiplication by $\rho$'', $L = -\Delta_A$. Then $T$ has a non-null kernel, the space of functions $V \to \C$ which are null out of $\{x_0\}$. Let $H_0 = (\ker T)^\perp$, i.e.\ the space of functions which vanish in $x_0$, so that $F_r|_{H_0} = 0$ for $r \leq 1$. Then 
\[
\left.\begin{array}{r}
\|\chi_{\{\rho < t\}} f\|_2\\
\|e^{-\rho/t} f\|_2
\end{array}\right\} \leq C_{\gamma,\delta} \,t^{\delta} \|(-\Delta_A)^{\gamma} f\|_2 \qquad \text{for $t$ small}
\]
trivially holds for $f \in H_0$.

We consider now two particular cases. The first one is the $n$-dimensional square lattice, with $V = \Z^n$ and edges only between vertices $(x_1,\dots,x_n)$, $(y_1,\dots,y_n)$ such that $\sum_{j=1}^n |x_j - y_j| = 1$. By direct calculation through Fourier series, one obtains
\begin{multline*}
\|E_r\|_{1 \to \infty} = \lambda^n \left(\left\{x \in \left[-\frac{1}{2},\frac{1}{2}\right]^n \tc \sum_{i=1}^n \left(1-\cos (2 \pi x_i) \right) < \frac{r}{2}\right\}\right)\\
\asymp \begin{cases}
r^{n/2} &\text{for $r \to 0^+$}\\
1 &\text{for $r \to +\infty$}
\end{cases}
\end{multline*}
(where $\lambda^n$ is Lebesgue measure in $\R^n$), whereas
\[\|F_r\|_{\infty \to 1} = m(B(x_0,r)) \asymp \begin{cases}
1 &\text{for $r \to 0^+$}\\
r^n &\text{for $r \to +\infty$}
\end{cases}\]
Therefore Theorem~\ref{thm:uncertainty1} can be applied with $L,T$ swapped, $\Phi(r) = c r^{n/2}$ on the interval $\left]0,1\right[$ to obtain: for $0 < \gamma < n/2$, $f \in L^2(G)$,
\[\left.\begin{array}{r}
\|\chi_{\{\rho < t\}} f\|_2\\
\|e^{-\rho/t} f\|_2
\end{array}\right\} \leq C_\gamma t^{\gamma} \|(-\Delta_A)^{\gamma/2} f\|_2 \qquad \text{for $t$ large.}\]
From this and Theorem~\ref{thm:uncertainty2}, restricted global inequalities follow:

\begin{cor}
If $G$ is the $n$-dimensional square lattice, for $\alpha,\beta > 0$, $f \in L^2(G)$ with $f(0) = 0$,
\[\|f\|_2 \leq C_{\alpha,\beta} \|\rho^\alpha f\|_2^\frac{\beta}{\alpha+\beta} \|(-\Delta_A)^{\beta/2} f\|_2^\frac{\alpha}{\alpha+\beta}.\]
\end{cor}

Note that this inequalities can also be obtained from the corresponding inequalities for tori $\T^n$, which are a particular case of compact Riemannian manifolds. In fact, through Fourier transform, $H_0$ on $\Z^n$ corresponds to the space of functions with null mean value on $\T^n$, multiplication by $-\rho^2$ on $\Z^n$ corresponds to the Laplacian on $\T^n$, $-\Delta_A$ on $\Z^n$ corresponds to multiplication by
\[2 \sum_{i=1}^n (1-\cos(2\pi x_i)) \asymp \sum_{i=1}^n x_i^2\]
on $\T^n$.

The second case which we consider is the homogeneous tree of degree $n$, with $n > 2$ (note that the tree with $n = 2$ coincides with the $1$-dimensional square lattice). In this case, the spectrum of the adjacency matrix $A$ is known to be $\left[-2\sqrt{n-1},2\sqrt{n-1}\right]$, so that (since $n > 2$) $L$ has a spectral gap ($E_r = 0$ for $r < b = n - 2\sqrt{n-1}$) and, as in the case of Riemannian manifolds, local inequalities for large times, but also restricted global inequalities become trivial (since $L$, $T|_{H_0}$ have bounded inverses). A more interesting result is obtained by replacing $L$ with $\tilde L = L - b$. In fact, it is known (see~\cite{MR1653343}) that
\[\|e^{-tL}\|_{1 \to \infty} \asymp t^{-3/2} e^{-bt} \qquad\text{for $t$ large},\]
whereas
\[m(B(x_0,r)) \asymp (n-1)^r = e^{\kappa r} \qquad\text{for $r$ large}\]
(where $\kappa = \log(n-1)$), so that
\[\|\tilde E_{1/t}\|_{1 \to \infty} \lesssim \|e^{2t\tilde L}\|_{1 \to \infty} \asymp t^{-3/2} \qquad\text{for $t$ large}\]
and, putting $\tilde\rho = e^{\frac{\kappa}{3}\rho}$,
\[m(\{\tilde\rho < r\}) \lesssim r^3 \qquad\text{for $r$ large}.\]
Therefore Theorem~\ref{thm:uncertainty1} can be applied to $\tilde\rho, \tilde L$, obtaining
\[\left.\begin{array}{r}
\|\chi_{\{\tilde\rho < t\}} f\|_2\\
\|e^{-\tilde\rho/t} f\|_2
\end{array}\right\} \leq C_\gamma t^{\gamma} \|\tilde L^{\gamma/2} f\|_2 \qquad \text{for $t$ large,}\]
for $\gamma < 3/2$, $f \in L^2(G)$. Since this inequality trivially holds for $t$ small, by Theorem~\ref{thm:uncertainty2} we get uncertainty inequalities for $\tilde\rho, \tilde L$:

\begin{cor}
If $G$ is the homogeneous tree of degree $n$, for all $\alpha,\beta > 0$, $f \in L^2(G)$,
\[\|f\|_2 \leq C_{\alpha,\beta} \|e^{\alpha \frac{\log(n-1)}{3} \rho} f\|_2^\frac{\beta}{\alpha+\beta} \|(L-b)^{\beta/2} f\|_2^\frac{\alpha}{\alpha+\beta}.\]
\end{cor}

\subsection{Unimodular Lie groups and sublaplacians}

Results about the Laplace-Beltrami operator can be generalized to sublaplacians. In order to obtain uniform estimates, we restrict here to the case of left-invariant sublaplacians on connected unimodular Lie groups (see~\cite{MR924464}, \cite{MR1218884} for a reference).

Let $G$ be a connected unimodular Lie group, $m$ a Haar measure, $H = \{X_1,\dots,X_k\}$ a system of left-invariant vector fields on $G$ satisfying the H\"ormander condition, $L = - \sum_{i=1}^k X_i^2$ the associated sublaplacian. Then $L$ is a positive symmetric linear operator on $C_c^\infty(G) \subseteq L^2(G)$ and we can consider its Friedrichs extension, also denoted by $L$, which is positive self-adjoint on $L^2(G)$; its exponential semigroup $e^{-tL}$ ($t > 0$) admits moreover a kernel function $h_t$, the heat kernel, which has the same properties listed in \S\ref{subsection:riemannian} for the Riemannian case.

Let $d$, $\delta$ be respectively the Carnot-Carath\'eodory distance and the local dimension associated to $H$, $x_0 \in G$, $\rho = d(x_0,\cdot)$, $T$ the operator ``multiplication by $\rho$''. Then, for $r,t > 0$ small,
\[\|F_r\|_{\infty \to 1} = m(B(x_0,r)) \asymp r^\delta, \qquad\|E_{1/t}\|_{1 \to \infty} \asymp \|h_{2t}\|_\infty \asymp t^{-\delta/2},\]
so that local uncertainty inequalities can be obtained as in the Riemannian case.

To extend such inequalities to large times, it is useful to recall a result of Guivarc'h \cite{MR0369608}, which states that the volume growth of $G$ can be either strictly polynomial:
\[m(B(x_0,r)) \asymp r^a \qquad\text{for some $a \in \N$ and for $r \to +\infty$}\]
or exponential:
\[e^{\beta r} \lesssim m(B(x_0,r)) \lesssim e^{\kappa r} \qquad\text{for some $\beta,\kappa > 0$ and for $r \to +\infty$}.\]

In the polynomial case, it is known that
\[\|h_t\|_\infty \asymp t^{-a/2} \qquad\text{for $t \to +\infty$},\]
therefore, exactly as in the Riemannian case, global uncertainty inequalities can be obtained (this is one of the results of \cite{MR2355602}):

\begin{cor}
If $G$ is a connected unimodular Lie group with polynomial growth, for all $\alpha,\beta > 0$, $f \in L^2(G)$,
\[\|f\|_2 \leq C_{\alpha,\beta} \|\rho^\alpha f\|_2^\frac{\beta}{\alpha+\beta} \|L^{\beta/2} f\|_2^\frac{\alpha}{\alpha+\beta}\]
(except for the compact case, in which we have to restrict to the functions $f$ with null mean value).
\end{cor}

In the exponential case, instead,
\begin{equation}\label{eq:expgroupheatdecay}
\|E_{1/t}\|_{1 \to \infty} \lesssim \|h_{2t}\|_{\infty} \lesssim e^{-ct^{1/3}} \qquad\text{for $t \to +\infty$}
\end{equation}
for some $c > 0$. Putting
\[\Phi(r) = \begin{cases}
r^\delta &\text{if $r \leq 1$}\\
e^{\kappa(r-1)} &\text{if $r \geq 1$}
\end{cases}\]
we have that $\Phi$ satisfies \eqref{eq:ip1} for $\gamma < n/2$ and moreover
\[\|F_r\|_{\infty \to 1} \lesssim \Phi(r), \qquad \|E_{1/t}\|_{1 \to \infty} \,\Phi(c \kappa^{-1} t^{1/3}) \lesssim 1\]
for all $r > 0$ and for $t$ large. Therefore, by Theorem~\ref{thm:uncertainty1}, for $\gamma < \delta/2$, $f \in L^2(G)$,
\[\left.\begin{array}{r}
\|E_{1/t} f\|_2\\
\|e^{-tL} f\|_2
\end{array}\right\} \leq C_\alpha t^{-\gamma/3} \|\rho^{\gamma} f\|_2 \qquad\text{for $t$ large.}\]
Unfortunately, this local inequality cannot be combined with the one for small times, since $t^{-\gamma/3} < t^{-\gamma/2}$ for $t$ small and $t^{-\gamma/2} < t^{-\gamma/3}$ for $t$ large.

To obtain a global inequality, we can slightly modify the operators $T,L$. For instance, if we replace the distance function $\rho$ with
\[\tilde\rho = \rho(1+\rho)^{1/2}\]
we easily get
\[m(\{\tilde \rho < r\}) \lesssim \begin{cases}
r^\delta &\text{for $r$ small}\\
e^{\kappa r^{2/3}} &\text{for $r$ large}
\end{cases}\]
so that, by Theorem~\ref{thm:uncertainty1}, the inequality
\[\left.\begin{array}{r}
\|E_{1/t} f\|_2\\
\|e^{-tL} f\|_2
\end{array}\right\} \leq C_\alpha t^{-\gamma/2} \|\tilde\rho^{\gamma} f\|_2\]
holds for all times (and $\gamma < \delta/2$); therefore we obtain the following global inequality:

\begin{cor}
If $G$ is a connected unimodular Lie group with exponential growth, for all $\alpha,\beta > 0$, $f \in L^2(G)$,
\[\|f\|_2 \leq C_{\alpha,\beta} \|\rho^\alpha (1+\rho)^{\alpha/2} f\|_2^\frac{\beta}{\alpha+\beta} \|L^{\beta/2} f\|_2^\frac{\alpha}{\alpha+\beta}.\]
\end{cor}

It should be remarked that the estimate~\eqref{eq:expgroupheatdecay} is not always optimal: if $L$ has a spectral gap (i.e.\ if $G$ is not amenable, cf.~\cite{MR1355612}), we have $E_{1/t} = 0$ for $t$ large and we can proceed as in the Riemannian case. However, there do exist unimodular Lie groups with exponential volume growth and without spectral gap (for an example, see~\cite{MR1415763}).

The work of Varopoulos \cite{MR1355612} (cf.\ also \cite{MR1435484}) allows us to obtain more precise results in the case of non-amenable groups. Let $b$ be the spectral gap of $L$ (i.e.\ $\|e^{-tL}\|_{2 \to 2} = e^{-tb}$) and $Q$ be the radical of $G$; then, if $Q$ has polynomial growth,
\[\|h_t\|_{\infty} \lesssim t^{-\nu/2} e^{-b t} \qquad\text{for $t \geq 1$}\]
for some $\nu \geq 0$, whereas, if $Q$ has exponential growth,
\[\|h_t\|_{\infty} \lesssim e^{-b t - c t^{1/3}} \qquad\text{for $t \geq 1$}\]
for some $c > 0$. This means that, putting $\tilde L = L - b$, we have $e^{-t\tilde L} = e^{bt} e^{-tL}$, so that
\[\|e^{-t\tilde L}\|_{1 \to \infty} \lesssim \begin{cases}
t^{-\delta/2} &\text{($t$ small)}\\
t^{-\nu/2} &\text{($Q$ polynomial, $t$ large)}\\
e^{-ct^{1/3}} &\text{($Q$ exponential, $t$ large)}
\end{cases}\]
and, replacing $L$ with $\tilde L$, we can proceed as before.

Namely, if $Q$ has exponential growth, then we get

\begin{cor}
If $G$ is a non-amenable connected unimodular Lie group whose radical has exponential growth, for all $\alpha,\beta > 0$, $f \in L^2(G)$,
\[\|f\|_2 \leq C_{\alpha,\beta} \|\rho^\alpha (1+\rho)^{\alpha/2} f\|_2^\frac{\beta}{\alpha+\beta} \|(L-b)^{\beta/2} f\|_2^\frac{\alpha}{\alpha+\beta}.\]
\end{cor}

If on the contrary $Q$ has polynomial growth and $\nu > 0$, then we can replace the distance function $\rho$ with
\[\tilde\rho = e^{\frac{\kappa}{\nu} \rho} - 1\]
so that
\[m(\{ \tilde\rho < r \}) \lesssim \begin{cases}
r^\delta &\text{for $r$ small}\\
r^\nu &\text{for $r$ large}
\end{cases}\]
and finally

\begin{cor}
If $G$ is a non-amenable connected unimodular Lie group whose radical is noncompact and has polynomial growth, for all $\alpha,\beta > 0$, $f \in L^2(G)$,
\[\|f\|_2 \leq C_{\alpha,\beta} \|(e^{\frac{\kappa}{\nu}\rho} -1)^\alpha f\|_2^\frac{\beta}{\alpha+\beta} \|(L-b)^{\beta/2} f\|_2^\frac{\alpha}{\alpha+\beta}.\]
\end{cor}

\bibliographystyle{abbrv}
\bibliography{../licenza}

\end{document}